\documentclass[12pt,a4paper,reqno]{article}
\usepackage[scaled]{helvet}

\usepackage{graphicx}
\usepackage{epstopdf} 
\usepackage{amsfonts}
\usepackage{amsthm}
\usepackage{amsmath}
\usepackage{amscd}
\usepackage[mathscr]{eucal}
\usepackage{indentfirst}
\usepackage{pict2e}
\usepackage{epic}
\numberwithin{equation}{section}
\usepackage{MnSymbol}
\usepackage[top=1.5cm, bottom=1cm, left=1.5 cm, right=1cm]{geometry}
\usepackage{tikz}
\usepackage{enumitem}

\usepackage{hyperref}

\numberwithin{equation}{section}

\theoremstyle{plain}
\newtheorem{Th}{Theorem}[section]
\newtheorem{Lemma}[Th]{Lemma}

\newtheorem{Prop}[Th]{Proposition}

 \theoremstyle{definition}
\newtheorem{Def}[Th]{Definition}

\newtheorem{Rem}[Th]{Remark}
\newtheorem{?}[Th]{Problem}
\newtheorem{Ex}[Th]{Example}

\begin{document}
\vskip .2cm 
\begin{center}
\huge\textbf{Generalized Uncertainty Principles associated with the Quaternionic Offset Linear Canonical Transform}
\end{center}

\begin{center}
\textbf{Youssef El Haoui$^{1,}$\footnote[1]{Corresponding author.}, Said Fahlaoui$^1$, Eckhard Hitzer$^2$}

 $^1$Department of Mathematics and Computer Sciences, Faculty of Sciences, University Moulay Ismail, Meknes 11201, Morocco\\
$^2$Dr. rer. nat. of Theoretical Physics, International Christian University,\\ Osawa 3-10-2, Mitaka-shi 181-8585 Tokyo, Japan

E-MAIL: youssefelhaoui@gmail.com, saidfahlaoui@gmail.com, hitzer@icu.ac.jp 
\end{center}

\vspace{1cm}
\begin{abstract}The quaternionic  offset linear canonical transform (QOLCT) can be defined as a generalization of the quaternionic linear canonical transform (QLCT). 
In this paper, we define  the QOLCT, we derive the relationship between the QOLCT and the quaternion Fourier transform (QFT).
 Based on this fact, we prove the Plancherel formula, and some properties related to the QOLCT. Then, we generalize some different uncertainty principles (UPs), including Heisenberg-Weyl’s UP, Hardy’s UP, Beurling’s UP, and  logarithmic UP to the QOLCT domain in a broader sense.

\end{abstract}




\section{Introduction}
 The QFT plays an important role  in the representation of signals. It transforms a real (or quaternionic) 2D signal into a quaternion valued frequency domain signal.
In \cite{EF17}, the authors provide the Heisenberg's inequality and the Hardy's UP  for the two-sided QFT. The authors in \cite{FE17} generalize the Beurling's UP to the QFT domain.
It is well known that the LCT provides a more general framework for a number of famous linear integral transforms in signal processing and optics, such as Fourier transform FT, the fractional FT, the Fresnel transform, the Lorentz transform.\\
The LCT was extended to the Clifford analysis by  Kit Ian Kou et al \cite{KMZ13} in the 2013s, to study the generalized prolate spheroidal wave functions and their connection with energy concentration problems.\\
In \cite{KOM13}, the authors introduced the quaternion linear canonical transform (QLCT), which is a generalization of the LCT in the framework of quaternion algebra. \\
Several properties, such as the Parseval's formula, and UP associated with  the QLCT are established.\\
In view of the fact that the OLCT is a generalization of the LCT, and  has wide applications in signal processing and optics, one is interested to extend the  OLCT  to a quaternionc algebra framework.\\
To the best of our knowledge, the generalization of the OLCT to a quaternionic  algebra, and the study of the properties and UPs associated with this generalization, have not been carried out yet.
 Therefore, the results in this paper are new in the literature.\\
 The main objective of the present study is to develop further technical methods in the theory of partial differential equations \cite{FE83}.\\
In the present work, we study the QOLCT that transforms a real (or quaternionic) 2D signal into a quaternion-valued frequency domain signal. Some important properties of the two-sided QOLCT are established.  A well known UPs for the two-sided QOLCT  are generalized.

The rest of the paper is organized as follows: Section 2 gives a brief introduction to some general definitions and basic properties of quaternionic analysis, and contains a reminder of the definition and some results for the two-sided QFT useful in the sequel.\\
The QOLCT of 2D quaternionic signal is introduced and studied in Section 3. Some important properties such as Plancherl’s  theorem are obtained, we also give the QOLCT of a Gaussian quaternionic functions (Gabor filters) to be indeed the only functions that minimize  the Heisenberg-Weyl's UP  associated with  the QOLCT, which has been proven
in section 4. In this section, we generalize the corresponding results of Hardy’s UP, Beurling’s UP, and  logarithmic UP to the QOLCT domain respectively. 
In section 5, we conclude the paper.

%

\section{Preliminaries}
\subsection*{The quaternion algebra} \ \\
In the present section we collect some basic facts about quaternions, which will be needed throughout the paper.
For all what follows, let $\mathbb{H}$ be the Hamiltonian skew field of quaternions:
$\mathbb{H}=\{q=q_0+iq_1+jq_2+kq_3;\ q_0, q_1, q_2, q_3 \in \mathbb{R}$\}\\
which is an associative noncommutative four-dimensional algebra.\\
where the elements ${  i},{  \ }{  j},{  \ }{  k}$ satisfy the Hamilton's multiplication rules:

$ij = -ji = k;\ jk = -kj = i;\ ki = -ik = j;\ i^2 = j^2 = k^2=-1.$
In this way the quaternionic algebra can be seen as an extension of the complex field $\mathbb{C}$.

Quaternions are isomorphic to the Clifford algebra ${Cl}_{(0,2)}$ of ${\mathbb R}^{(0,2)}$:

\begin{equation}\label{isom}
\mathbb{H} \cong {Cl}_{(0,2)}
\end{equation}

The scalar part of  a quaternion  $q \in \mathbb{H}$\ is $q_0$ denoted by $Sc(q)$, the non scalar part(or pure quaternion) of $q$\ is $iq_1+jq_2+kq_3$ denoted by $Vec(q)$.

The quaternion conjugate of $q \in \mathbb{H}$, given by

$\overline{q}=q_0-iq_1-jq_2-kq_3$. \\ is an anti-involution, namely,

\[\overline{qp}= \overline{p} \ \overline{q},\ \overline{p+q}= \overline{p}+\overline{q},\ \overline{\overline{p}}=p.\] 

The norm or modulus of $q \in {\mathbb H}$\ is defined by

\[{|q|}_Q=\sqrt{q\overline{q}}=\sqrt{{q_0}^2+\ {q_1}^2+{q_2}^2+{q_3}^2}.\] 

Then, we have

\[{|pq|}_Q={|p|}_Q{|q|}_Q.\] 
In particular, when $q=q_0$ is a real number, the module ${|q|}_Q$ reduces to the ordinary Euclidean module $\left|q\right|=\sqrt{{q_0}^2}$.

It is easy to verify that  $0 \ne q \in \mathbb{H}$ implies~:

\[q^{-1}=\frac{\overline{q}}{{{|q|}_Q}^2}.\] 

Any quaternion  $q$ can be written as $q$=\ ${{|q|}_Qe}^{\mu \theta }$  where $e^{\mu \theta }$ is understood in accordance with Euler's formula\\
 $e^{\mu \theta }={\cos  \left(\theta \right)\ }+\mu \ {\sin  \left(\theta \right),\ }$ where 
$\theta = artan \frac{\left|Vec\left(q\right)\right|_Q}{Sc\left(q\right)}$,\ 0$\le \theta \le \pi $ and \ $\mu $ :=\ $\frac{Vec\left(q\right)}{\left|{  Vec}\left({  q}\right)\right|_Q}$ verifying ${\mu }^2 =\ -1$.

In this paper, we will study the  quaternion-valued signal $f:{\mathbb{R}}^2\to \mathbb{H} $, $f$ which can be expressed as $$f=f_0+i f_1+jf_2+kf_3,$$ with $f_m~:  {\mathbb R}^2 \to \ {\mathbb R}\ for\ m=0,1,2,3.$
Let us introduce the canonical inner product for quaternion valued functions  $f,g\ :{\mathbb R}^2\ \to {\mathbb H}$, as follows:
\[<f,g> = \int_{{\mathbb R}^2} {f\left(t\right)\overline{g\left(t\right)}}dt,\   dt={dt}_1{dt}_2.\] 
Hence, the natural norm is given  by
\[{\left|f\right|}_{2,Q}=\sqrt{<f,f>} ={(\int_{{\mathbb R}^2}{{\left|f(t)\right|}^2_Q}dt)}^{\frac{1}{2}},\] 

and the quaternion module $L^2({\mathbb R}^2,\ {\mathbb H})$, is given by

\[L^2({\mathbb R}^2,\ {\mathbb H}) = \{f : {\mathbb {\mathbb R}^2\ \to {\mathbb H},\ {\left|f\right|}_{2,Q}< \infty }\}.\]

We denote by ${\mathcal S}(\mathbb R^2,\mathbb H)$, the quaternion Schwartz space of $C^{\infty }$- functions $f$, from ${{ \mathbb R}}^2$ to $\mathbb H$, such that for all $m,n \in \mathbb N$\\
 \[{sup}_{t\in \mathbb R^2,{{\alpha}_1+{\alpha}_2\le n}}
 {({\left(1+\left|t\right|\right)}^m{\left|\frac{{\partial }^{{\alpha }_1+{\alpha }_2}}{{{\partial t}_1}^{{\alpha }_1}{{\partial t}_2}^{{\alpha }_2}}f(t)\right|}_Q)}<\infty, where \ ({\alpha }_1,{\alpha }_1)\in \mathbb N^2.\]
Besides the quaternion units $i, j, k$, we will use the following real vector notation:\\
$t=(t_1,t_2) \in {\mathbb R}^2, \ |t|^2= {t_1}^2+{t_2}^2, \ f(t) = f(t_1,t_2),\ dt={dt}_1{dt}_2,$\ and so on.\\

\subsection{The general two-sided QFT } 
The QFT which has been defined by Ell \cite{EL93}, is a generalization of the classical Fourier transform (CFT) using a quaternionic algebra framework.
Several known and useful properties, and theorems of this extended transform are
generalizations of the corresponding properties, and theorems of the CFT with some
modifications (e.g., \cite{BU99}, \cite{CKL15}, \cite{HI07}, \cite{EL93}).
The QFT belongs to the family of Clifford Fourier transformations because of \eqref{isom}.
There are three different types of QFT, the left-sided QFT , the right-sided QFT, and two-sided QFT \cite{PDC01}.

Let us define the two-sided QFT and provide some properties used in the sequel. 
\begin{Def}
(Two-sided QFT with respect to two pure unit quaternions $\lambda;\mu$ \cite{HS13})\\
Let\ $\lambda, \mu \in  {\mathbb H}, {\lambda }^2= {\mu }^2=-1$, be any two pure unit quaternions.\\
For  $f$ in    $L^1\left({\mathbb R}^2,{\mathbb H}\right)$,\ the two-sided QFT with respect to $\lambda ; \mu $ is\\

\begin{equation}\label{QFT}
{\mathcal F}^{\lambda, \mu }\{f\}(u) =\int_{{\mathbb R}^2}{e^{-\lambda {{ u}}_1t_1}}\ f(t)\ e^{-\mu {{ u}}_2t_2}dt,~~ where ~t, u \in {\mathbb R}^2.
\end{equation}

\end{Def}

We define a new module of $\mathcal F\{f\}^{\lambda, \mu } $ \ as follows :
\begin{equation}\label{norm1}
{\left\|{\mathcal F}^{\lambda, \mu }\left\{f\right\}\right\|}_Q\ := \sqrt{\sum^{m=3}_{m=0}{{\left|{\mathcal F}^{\lambda, \mu }\left\{f_m\right\}\right|}^2_Q}}.\end{equation}

Furthermore, we define a new $L^2$-norm of ${\mathcal F\{f\}}$ as follows :\\
\begin{equation}\label{norm2}
{\left\|{\mathcal F}^{\lambda, \mu }\{f\}\right\|}_{2,Q}:=\sqrt{\int_{\mathbb {\mathbb R}^2}{{\left\|{\mathcal F}^{\lambda, \mu }\left\{f\right\}(y)\right\|}^2_Qdy}}.\end{equation}
It is interesting to observe that ${\left\|{\mathcal  F}^{\lambda,\mu }\{f\}\right\|}_{Q}$    is not equivalent to$\ {\left|{\mathcal  F}^{\lambda,\mu }\{f\}\right|}_{Q}$  unless  $f$ is real valued.\\

\begin{Lemma}\label{dilation}{(Dilation property), see example 2 on page 50 \cite{BU99}} \ \\
Let $k_1, k_2$ be a positive scalar constants,  we have

\[{\mathcal F}^{\lambda,\mu}\left\{f(t_1,t_2)\right\}\left(\frac{u_1}{k_1},\frac{u_2}{k_1}\right)=k_1k_2{\mathcal F}^{\lambda,\mu}\left\{f(k_1t_1,k_2t_2)\right\}\left(u_1,u_2\right).\] 
\end{Lemma}

By following the proof of theorem 3.2 in \cite{CKL15}, and replacing $i$ by $\lambda$,\ $j$ by $\mu$\ we obtain the next lemma.

\begin{Lemma}\label{Plancherel}{(QFT Plancherel)}\\
Let $f \in L^2({\mathbb R}^2,{\mathbb H})$, then 

\[\int_{{\mathbb R}^2}{{\left\|{\mathcal F}^{\lambda,\mu }\left\{f\right\}\left(u\right)\right\|}^2_Q}du=4{\pi }^2\int_{{\mathbb R}^2}{{\left|f(t)\right|}^2_Q}dt.\] 

\end{Lemma}

\begin{Lemma}\label{derivation} \ \\
 If   $f\in L^2({\mathbb R}^2,{\mathbb H}),\ \frac{{\partial }^{m+n}}{{\partial t}^m_1{\partial t}^n_2}f$ exist and are in   $L^2\left({\mathbb R}^2,{\mathbb H}\right)$ for \ $m,n\in {\mathbb N}_0,$\  
 then
\[{\mathcal F}^{\lambda,\mu }\left\{\frac{{\partial }^{m+n}}{{\partial t}^m_1{\partial t}^n_2}f\right\}\left(u\right)={(\lambda u_1)}^m\ {\mathcal F}^{\lambda,\mu }\left\{f\right\}\left(u\right)  {(\mu u_2)}^n.\] 
\end{Lemma}
Proof. See [\cite{BU99}, Thm. 2.10].

\begin{Lemma}\label{inverse}{Inverse QFT} \cite{HI14}

If  $f\in L^1\left({\mathbb{R}}^2, {\mathbb{H}}\right), and \ {\mathcal F}^{\lambda,\mu}\{f\}\in L^1\left({\mathbb{R}}^2,{\mathbb{H}}\right)$, then the two-sided QFT is an invertible transform and its inverse is given
by\\
\[f(t)= \frac{1}{{(2\pi )}^2} \int_{{\mathbb{R}}^2}{e^{\lambda u_1 t_1}} {\mathcal F}^{\lambda,\mu} \{f(t)\}(u) e^{\mu u_2 t_2}du. \]

\end{Lemma}
\section{The offset quaternionic linear canonical transform }

Kit Ian Kou et al \cite{KMZ13} introduce the quaternionic linear canonical transform (QLCT). They consider a pair of unit determinant two-by-two matrices 
\[A_1=\left[ \begin{array}{cc}
a_1 & b_1 \\ 
c_1 & d_1 \end{array}
\right], A_1=\left[ \begin{array}{cc}
a_2 & b_2 \\ 
c_2 & d_2 \end{array}
\right]\in  {\mathbb R}^{2\times 2},\] 

with unit determinant, that is \ $a_1d_1-b_1c_1=1,\ a_2d_2-b_2c_2=1,$\\
Eckhard Hitzer \cite{HI2014} generalize the definitions of \cite{KMZ13} to be:\\ 
The two-sided QLCT of signals f $\in L^1({\mathbb R}^2,{\mathbb H})$, is defined as\\
\[{\mathcal L}^{\lambda,\mu }_{A_1,A_2}\{f\}(u)= \int_{{\mathbb R}^2}{K^{\lambda }_{A_1}\left(t_1,u_1\right)}f\left(t\right)K^{\mu}_{A_2}\left(t_2,u_2\right)dt.\]
with $\lambda,\mu \in {\mathbb H},$\ two pure unit quaternions, ${\lambda }^2={\mu }^2=-1$,\ including the cases $\lambda =\pm \mu,$\\
$K^{\lambda }_{A_1}\left(t_1,u_1\right)=\frac{1}{\sqrt{\lambda 2\pi b_1}} e^{{\lambda (a_1t^2_1-2t_1u_1+d_1u^2_1)}/{2b_1}}, \  \ K^{\mu }_{A_2}\left(t_2,u_2\right)=\frac{1}{\sqrt{\mu 2\pi b_2}} e^{{\mu (a_2t^2_2-2t_2u_2+d_2u^2_2)}/{2b_2}},$\\
In \cite{KMZ13}, for $\lambda  = i$\ and\  $\mu = j$, the right-sided QLCT and its properties, including an UP are studied in some detail.

We now generalize the definitions of \cite{HS13}, \cite{HI14} as follows:
\begin{Def}
Let  $A_l=\left[\left| \begin{array}{cc}
a_l & b_l \\ 
c_l & d_l \end{array}
\right| \begin{array}{c}
{\tau }_l \\ 
{\eta }_l \end{array}
\right]$, 

parameters  $a_l, b_l, c_l,d_l,\ {\tau }_l,\ {\eta }_l\in {\mathbb R}$ such as $a_ld_l-b_lc_l=1$, for $l=1,2.$

The two-sided quaternionic offset linear canonical transform (QOLCT) of a signal $f \in L^1({\mathbb R}^2, {\mathbb H})$, is given by 

\begin{equation}\label{qolct}
{\mathcal O}^{\lambda,\mu }_{A_1,A_2}\{f(t)\}(u)=\left\{
 \begin{array}{c}

 \int_{{\mathbb R}^2}{K^{\lambda }_{A_1}\left(t_1,u_1\right)}f\left(t\right)K^{\mu }_{A_2}\left(t_2,u_2\right)dt,\ \ \ \ \  \ \ \ \ \ \ \  \ \  \ \ \ \  \  \ \ \ \ \ \ \ \   \   b_1,b_1\ne 0, \\ 
\sqrt{d_1}e^{\lambda (\frac{c_1d_1}{2}{\left(u_1-\tau_1\right)}^2+{\ u}_1\tau_1)}f\left(d_1{(u}_1-\tau_1\right),t_2)K^{\mu }_{A_2}\left(t_2,u_2\right),\ \ \  b_1=0,b_2\ne 0, \\ 
\sqrt{d_2}{K^{\lambda }_{A_1}\left(t_1,u_1\right)}f({t_1,d}_2{(u}_2-\tau_2))e^{\mu (\frac{c_2d_2}{2}{(u_2-\tau_2)}^2+{\ u}_2{\tau }_{2})},\ \ \ \ \ \ \ \ \ \ \  \ \ \ \ \ \ b_1\ne 0, b_2=0, \\ 
\sqrt{d_1d_2}{f(\left(d_1(u_1-{\tau}_1\right),d}_2(u_2-\tau_2)) e^{\lambda(\frac{c_1d_1}{2}{\left(u_1-\tau_1\right)}^2+{u}_1\tau_1)}e^{\mu(\frac{c_2d_2}2{(u_2-\tau_2)}^2+{u}_2{\tau}_2)}, \ b_1=b_2=0, 
\end{array}
\right.
\end{equation}

  $K^{\lambda }_{A_1}\left(t_1,u_1\right)=\frac{1}{\sqrt{\lambda 2\pi b_1}}\ e^{\lambda (a_1t^2_1-2t_1(u_1-\tau_1)-2{u_1(d}_1\tau_1-b_1\eta_1)+d_1{(u}^2_1+{\tau }^2_1))\frac{1}{2b_1}}$,  for $b_1\ne 0$,\\
and\\
$K^{\mu }_{A_2}\left(t_2,u_2\right)=\ \frac{1}{\sqrt{\mu 2\pi b_2}}\  e^{\mu (a_2t^2_2-2t_2(u_2-\tau_2)-2{u_2(d}_2\tau_2-b_2\eta_2)+d_2{(u}^2_2+{\tau }^2_2))\frac{1}{2b_2}}$, \ for\ $b_2\ne 0,$ 

with $\frac{1}{\sqrt{\lambda }}=e^{-\lambda \frac{\pi }{4}}$, $\frac{1}{\sqrt{\mu }} =e^{-\mu \frac{\pi }{4}}.$

The left-sided and right-sided QOLCTs can be defined correspondingly by placing the two kernel factors both on the left or on the right, respectively.
\end{Def}
We note that when $\tau_1=\tau_2=\ \eta_1=\eta_2=$0, the two-sided QOLCT reduces to the QLCT.

Also, when    $A_1=A_2=\left[\left| \begin{array}{cc}
0 & 1 \\ 
-1 & 0 \end{array}
\right| \begin{array}{c}
0 \\ 
0 \end{array}
\right]$, the conventional two-sided QFT is recovered. Namely,

${\mathcal O}^{\lambda,\mu }_{A_1,A_2}\{f(t)\}(u) =\frac{1}{\sqrt{\lambda 2\pi }}(\int_{{\mathbb R}^2}{ e^{-\lambda t_1u_1}}f\left(t\right)e^{-\mu t_2u_2}dt)\frac{1}{\sqrt{\mu 2\pi }}$

\hspace*{3 cm}=$\frac{1}{2\pi }e^{-\lambda \frac{\pi }{4}}{\mathcal F}^{\lambda,\mu }\left\{f\right\}\left(u_1,u_2\right)\ e^{-\mu \frac{\pi }{4}},$

where ${\ \mathcal F}^{\lambda,\mu }\left\{f\right\}$ is the QFT of \ $f$\ given by \eqref{QFT}.

For simplicity's sake, in this paper we restrict our attention to the two-sided QLCTs of 2D quaternion-valued signals.
Note that when $b_1b_2=$\textit{ }0 or $b_1=b_2=$\textit{ }0  the QOLCT of a function is essentially a chirp multiplication and is of no particular interest in our objective interests.
Hence, we deal with only the case when $b_1b_2 \ne 0$\  in this paper, without loss of generality, we set $b_l>0 (l = 1, 2),$  

The following lemma gives the relationships of two-sided QOLCTs and two-sided QFTs of 2D quaternion-valued signals.

\begin{Lemma}\label{QFT_QOLCT}
The QOLCT of a signal $f \in L^1({\mathbb R}^2, {\mathbb H})$\ can be reduced to the QFT

\[{\mathcal O}^{\lambda,\mu }_{A_1,A_2}\left\{f\left(t\right)\right\}\left(u_1,u_2\right)= {\mathcal F}^{\lambda,\mu}\left\{h(t)\right\}\left(\frac{u_1}{b_1},\frac{u_2}{b_2}\right),\ \] 

with \\
\begin{equation}\label{function_h}
h(t)= \frac{1}{\sqrt{2\pi {\lambda b}_1}}e^{\lambda[-\frac{1}{b_1}{u_1(d}_1\tau_1-b_1\eta_1)+\frac{d_1}{2b_1}{(u}^2_1+{\tau }^2_1)+\frac{1}{b_1}t_1\tau_1+\frac{a_1}{2b_1}t^2_1]}f(t)e^{\mu[-\frac{1}{b_2}{u_2(d}_2\tau_2-b_2\eta_2)+\frac{d_2}{2b_2}{(u}^2_2+{\tau }^2_2)+ \frac{1}{b_2}t_2{\tau}_2+\frac{a_2}{2b_2}t^2_2]}\frac{1}{\sqrt{2\pi \mu b_2}}.
\end{equation}
\end{Lemma}
Proof.
From the definition of the QOLCT, we have 

${\mathcal O}^{\lambda,\mu}_{A_1,A_2}\left\{f\left(t\right)\right\}\left(u_1,u_2\right)=\int_{{\mathbb R}^2}{K^{\lambda }_{A_1}\left(t_1,u_1\right)}f\left(t\right)K^{\mu }_{A_2}\left(t_2,u_2\right)dt$

=$\int_{{\mathbb R}^2}{\frac{1}{\sqrt{2\pi{\lambda b}_1}}e^{\lambda[\frac{a_1}{2b_1}t^2_1-\frac{1}{b_1}t_1(u_1-\tau_1)-\frac{1}{b_1}{u_1(d}_1\tau_1-b_1\eta_1)+\frac{d_1}{2b_1}{(u}^2_1+{\tau }^2_1)]}}f(t)\frac{1}{\sqrt{2\pi {\mu b}_2}}$\\
  \hspace*{2cm}        $ \times e^{\mu[\frac{a_2}{2b_2}t^2_2-\frac{1}{b_2}t_2(u_2-\tau_2)-\frac{1}{b_2}{u_2(d}_2\tau_2-b_2{\eta}_2)+\frac{d_2}{2b_2}{(u}^2_2+{\tau }^2_2)]}dt$

=$\int_{{\mathbb R}^2}{e^{-\lambda\frac{1}b_1t_1u_1}}[\frac{1}{\sqrt{2\pi\lambda b_1}}e^{\lambda[-\frac{1}{b_1}{u_1(d}_1\tau_1-b_1\eta_1)+\frac{d_1}{2b_1}(u^2_1+{\tau }^2_1)+\frac{1}{b_1}t_1\tau_1+\frac{a_1}{2b_1}t^2_1]}f(t)e^{\mu[-\frac{1}{b_2}{u_2(d}_2\tau_2-b_2\eta_2)+\frac{d_2}{2b_2}(u^2_2+{\tau }^2_2)+\frac{1}{b_2}t_2\tau_2+\frac{a_2}{2b_2}t^2_2]}\\ 
\hspace*{1cm} \times \frac{1}{\sqrt{2\pi{\mu b}_2}]}]  {e^{-\mu\frac{1}{b_2}t_2u_2}} dt$
 
\hspace*{0.1 cm}=${\mathcal F}^{\lambda,\mu }\left\{h(t)\right\}\left(\frac{u_1}{b_1},\frac{u_2}{b_2}\right)\ \hfill \square$

Due to the lemma 3.2 and  proposition 3.1, theorem 3.1 in \cite{CKL15}, the following properties are easily shown.
\begin{Th}

Let $f\in L^1({\mathbb R}^2,{\mathbb H})$. Then its QOLCT satisfies:

\hspace*{0.6 cm} $ \tikz\draw[black,fill=black] (0,0) circle (.5ex);$\ The map f  ${\to O}^{\lambda,\mu }_{A_1,A_2}\left\{f\right\}$ \ is real linear. That is, for $\alpha,\beta \in {\mathbb R}$, we have \\
\hspace*{2 cm}$ {\mathcal O}^{\lambda,\mu }_{A_1,A_2}\left\{\alpha f+\beta g\right\}=\alpha \ {\mathcal O}^{\lambda, \mu }_{A_1,A_2}\left\{f\right\}+\beta\  {\mathcal O}^{\lambda,\mu }_{A_1,A_2}\left\{g\right\}$.\\
\hspace*{1.2 cm} $ \tikz\draw[black,fill=black] (0,0) circle (.5ex);\ {\mathop{\lim }_{\left|u\right|\to \infty } {\left\|\ {\mathcal O}^{\lambda, \mu }_{A_1,A_2}\left\{f\right\}{ (u)}\right\|}_Q=0}.$

 \hspace*{0.7 cm}$\tikz\draw[black,fill=black] (0,0) circle (.5ex); \ {\mathcal O}^{\lambda,\mu }_{A_1,A_2}\left\{f\right\}$ is uniformly continuous on ${\mathbb R}^2.$

\end{Th}

Following the proof of theorems 11 and 12 in \cite{BA2016},  and by straightforward computation we derive shift and modulation properties for the QOLCT

\begin{Th} 

Let $f\in L^1({\mathbb R}^2,{\mathbb H}),$ with $t, u \in{\mathbb R}^2$,\ constants $\xi  =(\xi_1,\xi_2),\ k=(k_1,k_2) \in {\mathbb R}^2.$

We have:\\
\hspace*{0.6 cm} $ \tikz\draw[black,fill=black] (0,0) circle (.5ex);$\ t-Shift property\\
$O^{\lambda ,\mu }_{A_1,A_2}\{f(t-k)\}(u)=$\\
           $ e^{\lambda[({2k_1u_1-a}_1k^2_1)b_1c_1-2{k_1a_1(d}_1\tau_1-b_1\eta_1)]\frac{1}{2b_1}}O^{\lambda ,\mu }_{A_1,A_2}\{f(t)\}(u_1-{k_1a}_1,u_2-{k_2a}_2) e^{\mu[({2k_2u_2-a}_2k^2_2)b_2c_2-2{k_2a_2(d}_2\tau_2-b_2\eta_2)]\frac{1}{2b_2}}.$\\
\hspace*{0.6 cm} $ \tikz\draw[black,fill=black] (0,0) circle (.5ex);$\ Modulation  property\\
$O^{\lambda ,\mu }_{A_1,A_2}\{ e^{\lambda t_1\xi_1} f(t) e^{\mu t_2\xi_2}\}(u)=$\\
          $ e^{-\lambda  \left[\frac{d_1}{2}{(b}_1{\xi }^2_1-2b_1u_1\right)+{\xi_1(d}_1\tau_1-b_1\eta_1)]}O^{\lambda,\mu}_{A_1,A_2}\{f(t\}(u_1-{b_1\xi}_1,u_2-b_2\xi_2)e^{-\mu \left[\frac{d_2}{2}(b_2{\xi }^2_2-2b_2u_2\right)+\xi_2(d_2\tau_2-b_2\eta_2)]}.$
\end{Th}

\begin{Th}
If $f$  \  and ${\mathcal O}^{\lambda,\mu}_{A_1,A_2}\left\{f\right\}$ \ are in $ \ L^1({\mathbb R}^2,{\mathbb H}),$
then the inverse transform of the QOLCT can be derived from that of the QFT.
\end{Th} 

Proof. Indeed, Let \begin{equation}\label{function_g}
g\left(t\right)=e^{\lambda\frac{1}{b_1}t_1\tau_1+\lambda\frac{a_1}{2b_1}t^2_1}\ f(t) e^{\mu\frac{1}{b_2}t_2\tau_2+\mu\frac{a_2}{2b_2}t^2_2}.
\end{equation} 
We have \\
${\mathcal O}^{\lambda,\mu}_{A_1,A_2}\left\{f(t)\right\}\left(u_1,u_2\right)=\frac{1}{\sqrt{2\pi \lambda b_1}}e^{-\lambda\frac{1}{b_1}{u_1(d}_1\tau_1-b_1\eta_1)+\lambda\frac{d_1}{2b_1}{(u}^2_1+{\tau }^2_1)} {\mathcal F}^{\lambda,\mu}\left\{g(t)\right\}\left(\frac{u_1}{b_1},\frac{u_2}{b_2}\right)$\\ 
\hspace*{5cm}$\times e^{-\mu\frac{1}{b_2}{u_2(d}_2\tau_2-b_2\eta_2)+\mu\frac{d_2}{2b_2}{(u}^2_2+{\tau }^2_2)}\frac{1}{\sqrt{2\pi {\mu b}_2}}$ \\ 
\hspace*{3.3 cm}=$ e^{-\lambda(a_1t^2_1-{2t}_1(u_1-\tau_1))\frac{1}{2b_1}}K^{\lambda }_{A_1}\left(t_1,u_1\right){\mathcal F}^{\lambda,\mu}\left\{g(t)\right\}\left(\frac{u_1}{b_1},\frac{u_2}{b_2}\right)K^{\mu }_{A_2}\left(t_2,u_2\right)e^{-\mu(a_2t^2_2-{2t}_2(u_2-{\tau}_2))\frac{1}{2b_2}}.$\\
As $K_{A_m}\left(t_m,u_m\right) \overline{K_{A_m}\left(t_m,u_m\right)}=\frac{1}{2\pi b_m},\  m=1,2.$\  
We easily obtain\\
${\mathcal F}^{\lambda,\mu}\left\{g(t)\right\}\left(\frac{u_1}{b_1},\frac{u_2}{b_2}\right)={(2\pi)}^2b_1b_2\overline{K^{\lambda }_{A_1}\left(t_1,u_1\right)}e^{\lambda(a_1t^2_1-{2t}_1(u_1-\tau_1))\frac{1}{2b_1}}{\mathcal O}^{\lambda,\mu}_{A_1,A_2}\left\{f(t)\right\}\left(u_1,u_2\right)\overline{K^{\mu }_{A_2}\left(t_2,u_2\right)}$
\\ \hspace*{3 cm} $\times e^{\mu(a_2t^2_2-{2t}_2(u_2-\tau_2))\frac{1}{2b_2}}.$ \\
From lemma 2.5, it follows that 

\hspace*{1 cm}$g\left(t\right)=\frac{1}{{(2\pi )}^2}\int_{{\mathbb R}^2}{e^{\lambda t_1u_1}}{{\mathcal F}^{\lambda,\mu}\left\{g(t)\right\}\left(u\right)e}^{\mu t_2u_2}du$\\
\hspace*{2.50 cm}$=b_1b_2\int_{{\mathbb R}^2}{\overline{K^{\lambda }_{A_1}\left(t_1,{b_1u}_1\right)}e^{\lambda(a_1t^2_1+{2t}_1\tau_1)\frac{1}{2b_1}}{\mathcal O}^{\lambda,\mu}_{A_1,A_2}\left\{f(t)\right\}\left({b_1u}_1,b_2u_2\right)\overline{K^{\mu }_{A_2}\left(t_2,{b_2u}_2\right)}}$\\ 
 $\hspace*{3cm}\times e^{\mu(a_2t^2_2+{2t}_2\tau_2)\frac{1}{2b_2}}du.$\\ 
Or, equivalently\\
$e^{\lambda\frac{1}{b_1}t_1\tau_1+\lambda\frac{a_1}{2b_1}t^2_1}\ f(t)e^{\mu\frac{1}{b_2}t_2\tau_2+\mu\frac{a_2}{2b_2}t^2_2}= b_1b_2\int_{{\mathbb R}^2}\overline{K^{\lambda }_{A_1}\left(t_1,{b_1u}_1\right)}e^{\lambda(a_1t^2_1+{2t}_1\tau_1)\frac{a_1}{2b_1}}{\mathcal O}^{\lambda,\mu}_{A_1,A_2}\left\{f(t)\right\}\left({b_1u}_1,b_2u_2\right)$
\\ \hspace*{3 cm} $\times\overline{K^{\mu }{A_2}\left(t_2,{b_2u}_2\right)}e^{\mu(a_2t^2_2+{2t}_2\tau_2)\frac{a_2}{2b_2}}du.$\\
It means that
\hspace*{2 cm}$ f(t) =b_1b_2\int_{{\mathbb R}^2}{\overline{K^{\lambda }_{A_1}\left(t_1,{b_1u}_1\right)}{\mathcal O}^{\lambda,\mu}_{A_1,A_2}\left\{f(t)\right\}\left({b_1u}_1,b_2u_2\right)\overline{K^{\mu}_{A_2}\left(t_2,{b_2u}_2\right)}}du$\\
\hspace*{5.2cm}=$\int_{{\mathbb R}^2}{\overline{K^{\lambda }_{A_1}\left(t_1,u_1\right)}{\mathcal O}^{\lambda,\mu}_{A_1,A_2}\left\{f(t)\right\}\left(u_1,u_2\right)\overline{K^{\mu }_{A_2}\left(t_2,u_2\right)}}du.$

which is the inverse transform of the QOLCT. This proves the theorem.$\hfill \square$

\begin{Th}\label{Plancherel_QOLCT}{(Plancherel's theorem of the QOLCT)}\\
Every 2D quaternion-valued signal $f\in L^2({\mathbb R}^2,{\mathbb H})$ and its QOLCT are related to the Plancherel identity in the following way:

\[{\left\|{\mathcal O}^{\lambda,\mu }_{A_1,A_2}\left\{f\right\}\right\|}_{Q,2}={\left|f\right|}_{Q,2}.\] 
\end{Th}

Proof. Let $h(t)$ be rewritten in the form of \eqref{function_h}.

 By the definition of the norm ${ \left\|.\right\|}_{Q,2}$ and lemma \ref{dilation}  and lemma \ref{QFT_QOLCT}, we have

${\left\|{\mathcal O}^{\lambda,\mu }_{A_1,A_2}\left\{f\right\}\right\|}^2_{Q,2}=\int_{{\mathbb R}^2}{{\left\|\ {\mathcal O}^{\lambda,\mu }_{A_1,A_2}\left\{f\right\}{ (u)}\right\|}^2_Q}du$ 

\hspace*{2.4 cm}=$\int_{{\mathbb R}^2}{{\left\|{\mathcal F}^{\lambda,\mu }\left\{h(t)\right\}\left(\frac{u_1}{b_1},\frac{u_2}{b_2}\right)\right\|}^2_Q}du$

\hspace*{2.4 cm}=$\int_{{\mathbb R}^2}{{\left\|b_1b_2{\mathcal F}^{\lambda,\mu }\left\{h(b_1t_1,b_2t_2)\right\}\left(u_1,u_2\right)\right\|}^2_Q}du$

\hspace*{2.4 cm}=$b^2_1b^2_2\int_{{\mathbb R}^2}{{\left\|\ {\mathcal F}^{\lambda,\mu }\left\{h(b_1t_1, b_2t_2)\right\}\left(u_1,u_2\right)\right\|}^2_Q}du.$

From lemma \ref{Plancherel}  we get  

\[\int_{{\mathbb R}^2}{{\left\|\ {\mathcal F}^{\lambda,\mu }\left\{h(b_1t_1,b_2t_2)\right\}\left(u_1,u_2\right)\right\|}^2_Q}du={ 4}{\pi }^2\int_{{\mathbb R}^2}{{\left|h(b_1t_1,b_2t_2)\right|}^2_Q}dt.\] 

Let $s_l=b_lt_l,\ $\ for $l=1,2,$\ we have

\hspace*{1 cm} $ \int_{{\mathbb R}^2}{{\left|h(b_1t_1,b_2t_2)\right|}^2_Q}dt =\frac{1}{b_1b_2} \int_{{\mathbb {\mathbb R}}^2}{{\left|h(s_1,s_2)\right|}^2_Q}ds$

\hspace*{ 4.7cm}=$ \frac{1}{{{ 4}{\pi }^2b}^2_1b^2_2} \int_{{\mathbb R}^2}{{\left|f(s_1,s_2)\right|}^2_Q}ds.$

The last statement follows from ${\left|h\left(t\right)\right|}_Q= \frac{1}{2\pi \sqrt{b_1b_2}}{\left|\ f(t)\ \right|}_Q$,

Therefore, we get

\[{\left\|{\mathcal O}^{\lambda,\mu}_{A_1,A_2}\left\{f\right\}\right\|}^2_{Q,2}= \int_{{\mathbb R}^2}{{\left|f(s)\right|}^2_Q}ds={\left|f\right|}^2_{Q,2}.\] 

This ends the proof. $\hfill \square$

\begin{Lemma}

 If   $f\in L^2({\mathbb R}^2,{\mathbb H}),\ \frac{{\partial }^{l+n}}{{\partial t}^l_1{\partial t}^n_2}f$ exist and are in   $L^2\left({\mathbb R}^2,{\mathbb H}\right)$ for \ $l,n\in {\mathbb N}_0,$\  
 then\\
1. $\int_{{\mathbb R}^2}{u^2_1}{\left\| {\mathcal O}^{\lambda,\mu }_{A_1,A_2}\left\{f\right\}\left({ u}\right)\right\|}^2_Qdu={b^2_1\int_{{\mathbb R}^2}{}\left|\lambda (\frac{a_1}{b_1}t_1+\frac{\tau_1}{b_1}) f(t)+\frac{{\partial }}{{\partial t}_1} f(t)\right|}^2_Qdt,$\\
2. $\int_{{\mathbb R}^2}{u^2_2}{\left\|{\mathcal O}^{\lambda, \mu }_{A_1,A_2}\left\{f\right\}\left({ u}\right)\right\|}^2_Qdu={b^2_2\int_{{\mathbb R}^2}{}\left|(\frac{a_2}{b_2}t_2+ \frac{\tau_2}{b_2})f(t)\mu + \frac{{\partial }}{{\partial t}_2} f(t)\right|}^2_Qdt.$ 
\end{Lemma}

Proof. Let $h(t)$ be rewritten in the form of  \eqref{function_h}.
For the first statement, using lemma \ref{derivation} shows that
\[{\mathcal  F}^{\lambda,\mu }\left\{\frac{{\partial }}{{\partial t}_1}f\right\}\left(\frac{u_1}{b_1},\frac{u_2}{b_2}\right)=\lambda \frac{u_1}{b_1}{\mathcal  F}^{\lambda,\mu }\left\{f\right\}\left(\frac{u_1}{b_1},\frac{u_2}{b_2}\right),\] 
Then, using lemma \ref{QFT_QOLCT}, \eqref{norm1}, lemma \ref{Plancherel}, and the above equality we get\\
$\int_{{\mathbb R}^2}{u^2_1}{\left\|{\mathcal O}^{\lambda,\mu }_{A_1,A_2}\left\{f\right\}\left({ u}\right)\right\|}^2_Qdu= \int_{{\mathbb R}^2}{{u^2_1\ \left\|{\mathcal F}^{\lambda,\mu }\left\{h\left(t\right)\right\}\left(\frac{u_1}{b_1},\frac{u_2}{b_2}\right)\right\|}^2_Q}du$\\
\hspace*{4.7cm}=$\int_{{\mathbb R}^2}{\ \sum^{m=3}_{m=0}{{\left|{u_1\mathcal F}^{\lambda,\mu }\left\{h_m\left(t\right)\right\}\left(\frac{u_1}{b_1},\frac{u_2}{b_2}\right)\right|}^2_Q}}du$\\
\hspace*{4.7cm}=$\int_{{\mathbb R}^2}{\ \sum^{m=3}_{m=0}{{\left|\frac{1}{\lambda }b_1{\mathcal  F}^{\lambda,\mu }\left\{\frac{{\partial }}{{\partial t}_1}h_m\right\}\left(\frac{u_1}{b_1},\frac{u_2}{b_2}\right)\right|}^2_Q}}du$\\
\hspace*{4.7cm}= $b^2_1\int_{{\mathbb R}^2}{{\left\| {\mathcal  F}^{\lambda,\mu }\left\{\frac{{\partial }}{{\partial t}_1}h\right\}\left(\frac{u_1}{b_1},\frac{u_2}{b_2}\right)\right\|}^2_Q}du\ $\\
\hspace*{4.7cm}=$b^3_1b_2\int_{{\mathbb R}^2}{{\left\|{\mathcal  F}^{\lambda,\mu }\left\{\frac{{\partial }}{{\partial t}_1}h\right\}\left(\frac{u_1}{b_1},\frac{u_2}{b_2}\right)\right\|}^2_Q}\frac{{du}_1}{b_1}\frac{{du}_2}{b_2}$\\
\hspace*{4.7cm}=${4{\pi }^2b^3_1b_2\int_{{\mathbb R}^2}\left|\frac{{\partial }}{{\partial t}_1}h(t)\right|}^2_Qdt.$\\
where the last equation is the consequence of using lemma \ref{Plancherel}.\\
Moreover, using\\  
${\left|\frac{{\partial }}{{\partial t}_1}h\left(t\right)\right|}_Q=\\ \hspace*{1 cm}|\frac{{\partial}}{{\partial t}_1}(\frac{1}{\sqrt{2\pi {\lambda b}_1}}e^{\lambda[-\frac{1}{b_1}{u_1(d}_1\tau_1-b_1\eta_1)+\frac{d_1}{2b_1}{(u}^2_1+{\tau }^2_1)+\frac{1}{b_1}t_1\tau_1+\frac{a_1}{2b_1}t^2_1]}f(t)e^{\mu[-\frac{1}{b_2}{u_2(d}_2\tau_2-b_2\eta_2)+\frac{d_2}{2b_2}{(u}^2_2+{\tau }^2_2)+ \frac{1}{b_2}t_2{\tau}_2+\frac{a_2}{2b_2}t^2_2]}\frac{1}{\sqrt{2\pi \mu b_2}}) |_Q$\\
\hspace*{1.5cm}=$\frac{1}{2\pi \sqrt{b_1b_2}}|e^{\lambda [-\frac{1}{b_1}{u_1(d}_1\tau_1-b_1\eta_1)+\frac{d_1}{2b_1}{(u}^2_1+{\tau }^2_1)+\frac{1}{b_1}t_1\tau_1+\frac{a_1}{2b_1}t^2_1]}[\lambda (\frac{a_1}{b_1}t_1+\ \frac{\tau_1}{b_1}$)$\ f(t)$+$\ \frac{{\partial }}{{\partial t}_1}\ f(t)]|_Q$\\
\hspace*{1.5cm}=$\frac{1}{2\pi \sqrt{b_1b_2}}|\lambda (\frac{a_1}{b_1}t_1+\ \frac{\tau_1}{b_1}$)$\ f(t)+\ \frac{{\partial }}{{\partial t}_1}\ f(t)\ |_Q.$\\
We further get
\[\int_{{\mathbb R}^2}{u^2_1}{\left\|{\mathcal O}^{\lambda,\mu }_{A_1,A_2}\left\{f\right\}\left({ u}\right)\right\|}^2_Qdu={b^2_1\int_{{\mathbb R}^2}{}\left|\lambda (\frac{a_1}{b_1}t_1+\ \frac{\tau_1}{b_1}{ )}\ f(t)+\ \frac{{\partial }}{{\partial t}_1}\ f(t)\right|}^2_Qdt.\]

To prove the statement 2, we argue in the same spirit as in the previous proof.\\
Applying lemma \ref{QFT_QOLCT}, \eqref{norm1}, lemma \ref{derivation} and lemma \ref{Plancherel}, we have

$\int_{{\mathbb R}^2}{u^2_2}{\left\| {\mathcal O}^{\lambda,\mu }_{A_1,A_2}\left\{f\right\}({ u})\right\|}^2_Qdu=  \int_{{\mathbb R}^2}{{u^2_2\ \left\|{\mathcal F}^{\lambda,\mu }\left\{h(t)\right\}\left(\frac{u_1}{b_1},\frac{u_2}{b_2}\right)\right\|}^2_Q}du$\\
\hspace*{5.2cm}=$\int_{{\mathbb R}^2}{\sum^{m=3}_{m=0}{{\left|{\mathcal F}^{\lambda,\mu }\left\{h_m\left(t\right)\right\}\left(\frac{u_1}{b_1},\frac{u_2}{b_2}\right){\mu u}_2\right|}^2_Q}}du$\\
\hspace*{5.2cm}=$\int_{{\mathbb R}^2}{\ \sum^{m=3}_{m=0}{{\left| b_2 {\mathcal  F}^{\lambda,\mu }\left\{\frac{{\partial }}{{\partial t}_2}h_m\right\}\left(\frac{u_1}{b_1},\frac{u_2}{b_2}\right) \frac{1}{\mu }\right|}^2_Q}}du$\\
\hspace*{5.2cm}= $b^2_2\int_{{\mathbb R}^2}{{\left\|{\mathcal  F}^{\lambda,\mu }\left\{\frac{{\partial }}{{\partial t}_2}h\right\}\left(\frac{u_1}{b_1},\frac{u_2}{b_2}\right)\right\|}^2_Q}du $\\
\hspace*{5.2cm}=$b^3_2b_1\int_{{\mathbb R}^2}{{\left\|{\mathcal  F}^{\lambda,\mu }\left\{\frac{{\partial }}{{\partial t}_2}h\right\}\left(\frac{u_1}{b_1},\frac{u_2}{b_2}\right)\right\|}^2_Q}\frac{{du}_1}{b_1}\frac{{du}_2}{b_2}$\\
\hspace*{5.2cm}=${{ 4}{\pi }^2b^3_1b_2\int_{{\mathbb R}^2}{}\left|\frac{{\partial }}{{\partial t}_2}h(t)\right|}^2_Qdt.$

Since  ${\left|\frac{{\partial }}{{\partial t}_2}h\left(t\right)\right|}_Q=\frac{1}{2\pi \sqrt{b_1b_2}}|f(t)\mu (\frac{a_2}{b_2}t_2$+$\ \frac{\tau_2}{b_2}$)$\ $+$\ \frac{{\partial }}{{\partial t}_2}\ f(t)\ |_Q,$

It follows that \\

$\int_{{\mathbb R}^2}{u^2_2}{\left\| {\mathcal O}^{\lambda,\mu }_{A_1,A_2}\left\{f\right\}\left({ u}\right)\right\|}^2_Qdu={b^2_2\int_{{\mathbb R}^2}{}\left|(\frac{a_2}{b_2}t_2+ \frac{\tau_2}{b_2}{ )}f(t)\mu +\ \frac{{\partial }}{{\partial t}_2}\ f(t)\right|}^2_Qdt. \hfill \square$

\begin{Ex}{(The QOLCT  of a Gaussian quaternionic  function)}\\
Consider a Gaussian quaternionic  function \[f\left(t\right)=\beta e^{-({{\alpha }_1t}^2_1+{{\alpha }_2t}^2_2)},\] 

where $\beta ={\beta }_1{\beta }_2,$\ and\ ${\beta }_1={\beta }_{11}+\lambda{\beta }_{12},\  {\beta }_2={\beta}_{21}+\mu{\beta}_{22},\ $ and\ $ {\beta }_{11},{\beta }_{12},{\beta }_{21},{\beta }_{22}\in {\mathbb R},$ \ and\ ${\alpha }_1, {\alpha }_2$\ are real positive constants.\\

The QOLCT of $f$\ is given by\\
${\mathcal O}^{\lambda,\mu}_{A_1,A_2}\left\{f\left(t\right)\right\}\left(u_1,u_2\right)=\ {\beta }_1[\int_{{\mathbb R}^2}{K^\lambda_{A_1}\left(t_1,u_1\right)}e^{-{{\alpha }_1t}^2_1}e^{-{ {\alpha }_2t}^2_2}K^\mu_{A_2}\left(t_2,u_2\right)dt] {\beta }_2,$\\
\hspace*{3.8cm}  =${\beta }_1\int_{\mathbb R}{K^{\lambda}_{A_1}\left(t_1,u_1\right)}e^{-{{\alpha }_1t}^2_1}\ dt_1\int_{\mathbb R}{K^{\mu}_{A_2}\left(t_2,u_2\right)e^{-{\ {\alpha }_2t}^2_2}}dt\ {\beta }_2.$\\
We have
$\int_{\mathbb R}{K^{\lambda}_{A_1}\left(t_1,u_1\right)}e^{-{{\alpha }_1t}^2_1}\ dt_1=\frac{1}{\sqrt{2\pi\lambda b_1}}[\int_{\mathbb R}{\ e^{\lambda(\frac{a_1}{2b_1}t^2_1-t_1\frac{(u_1-\tau_1)}{b_1})}}e^{-{{\alpha }_1t}^2_1} dt_1]\ e^{\lambda(-{{2u}_1(d}_1\tau_1-b_1\eta_1)+d_1{(u}^2_1+{\tau }^2_1))\frac{1}{2b_1}},$\\
Since
 $\int_{\mathbb R}{ e^{\lambda\left(\frac{a_1}{2b_1}t^2_1-t_1\frac{\left(u_1-\tau_1\right)}{b_1}\right)}}e^{-{{\alpha }_1t}^2_1}\ dt_1=\int_{\mathbb R}{ e^{- \left({\alpha }_1- \frac{a_1}{2b_1}\lambda\right){[t_1+\lambda\frac{\frac{(u_1-\tau_1)}{b_1}}{2\left({\alpha }_1-\frac{a_1}{2b_1}\lambda\right)}]^2}}}\ dt_1  e^{-\frac{{(\frac{\left(u_1-\tau_1\right)}{b_1})}^2}{4\left({\alpha}_1-\frac{a_1}{2b_1}\lambda\right)}}$ \\
\hspace*{6.7 cm}=$\sqrt{\frac{\pi }{{\alpha }_1-\ \frac{a_1}{2b_1}\lambda}}e^{-\frac{{\left(u_1-\tau_1\right)}^2}{{2b}_1({2\alpha }_1b_1-a_1\lambda)}}{\alpha }^2_1$ \\    
\hspace*{6.7cm}=$\sqrt{\frac{2b_1\pi }{2b_1{\alpha }_1-a_1\lambda}} e^{-\frac{{\left(u_1-\tau_1\right)}^2}{{2b}_1\left({4\alpha }^2_1b^2_1+a^2_1\right)}({2\alpha }_1b_1+a_1\lambda)}$\\
\hspace*{6.7cm}=$\sqrt{\frac{2b_1\pi }{2b_1{\alpha }_1-a_1\lambda}} e^{- \frac{{\alpha }_1{\left(u_1-\tau_1\right)}^2}{\left({4\alpha }^2_1b^2_1+a^2_1\right)}}e^{-\frac{{a_1\left(u_1-\tau_1\right)}^2}{{2b}_1\left({4\alpha }^2_1b^2_1+a^2_1\right)}\lambda}$\\
where~the second equality follows from  $\int_{\mathbb R}{e^{-z{(t+z')}^2{\rm \ }}dt=}\sqrt{\frac{\pi }{z}}$,  for  $z,z'\in {\mathbb  C}$,Re(z)$>0,$\ (Gaussian integral with complex offset). \\
Then
$\int_{\mathbb R}{K^{\lambda }_{A_1}\left(t_1,u_1\right)}e^{-{{\alpha }_1t}^2_1}\ dt_1=\frac{1}{\sqrt{2b_1{\alpha }_1\lambda+a_1\lambda}} e^{-\frac{{\alpha }_1{\left(u_1-\tau_1\right)}^2}{\left({4\alpha }^2_1b^2_1+a^2_1\right)}}e^{\lambda(-{{2u}_1(d}_1\tau_1-b_1\eta_1)+d_1{(u}^2_1+{\tau }^2_1)-\frac{{a_1\left(u_1-\tau_1\right)}^2}{{4\alpha }^2_1b^2_1+a^2_1})\frac{1}{2b_1}}.$\\
We deduce that\\
${\mathcal O}^{\lambda,\mu}_{A_1\,A_1}\left\{f\left(t\right)\right\}\left(u_1{\rm,\ }u_2\right)= e^{- [\frac{{\alpha }_1{\left(u_1-\tau_1\right)}^2}{\left({4\alpha }^2_1b^2_1+a^2_1\right)}+\frac{{\alpha }_2{\left(u_2-\tau_2\right)}^2}{\left({4\alpha }^2_2b^2_2+a^2_2\right)}]}\ {\beta }_1\frac{1}{\sqrt{2b_1{\alpha }_1\lambda+a_1\lambda}}  
e^{\lambda(-{{2u}_1(d}_1\tau_1-b_1\eta_1)+d_1{(u}^2_1+{\tau }^2_1)-\ \frac{{a_1\left(u_1-\tau_1\right)}^2}{{4\alpha }^2_1b^2_1+a^2_1})\frac{1}{2b_1}} $ \\
\hspace*{7 cm}$\times\frac{1}{\sqrt{2b_2{\alpha }_2\mu+a_2\mu}}e^{\mu(-{{2u}_2(d}_2\tau_2-b_2\eta_2)+d_2{(u}^2_2+{\tau }^2_2)-\ \frac{{a_2\left(u_2-\tau_2\right)}^2}{{4\alpha }^2_2b^2_2+a^2_2})\frac{1}{2b_2}}{\beta }_2.\hfill \square$

\end{Ex}

Some properties of the QOLCT are summarized in Table \ref{QOLCT_prop}. 

\section{Uncertainty principles for the offset quaternionic linear canonical transform }
In harmonic analysis, the UP states that a non-trivial function and its  FT cannot both be sharply localized.  The UP plays an important role in signal processing, and quantum mechanics. 
In quantum mechanics , UP asserts that one cannot make certain of the position and momentum of  the particule at the same time, i.e., increasing the knowledge of the position decreases the knowledge of the momentum, and vice versa.
There are many different forms of UPs in the time-frequency plane, such as Heisenberg-Weyl's UP, Hardy's UP, Beurling's UP, and  logarithmic UP, and so on in terms of different notations of “localization”.
As far as we know, in  2013, Kit-Ian Kou et al \cite{KOM13} extended the Heisenberg-type UP to the QLCT. Recently Mawardi et al \cite{BA16} established  the logarithmic UP associated with the QLCT.
Considering that the QOLCT is a generalized version of the QLCT quaternionic Fourier,  and so of the QFT, it is natural and interesting to study
the simultaneous localization of a function and its QOLCT by further extending the aforementioned UPs to the QOLCT domain.
 Therefore, in this section, we prove and generalize the Heisenberg-Weyl's UP, Hardy's UP, Beurling's UP, and  logarithmic UP to 2D quaternion-valued signals using the two-sided QOLCT.

\subsection{Heisenberg-Weyl's uncertainty principle}

\begin{Prop}\label{QFT_Heisenberg}({\cite{EF17}, Thm. 4.1})\\

Let $f\left(t\right)={|f\left(t\right)|}_Qe^{u\left(t\right)\theta (t)}\ $.If $f,\ \frac{\partial }{\partial t_k} f,\ t_k f\in  L^2$($\mathbb R{^2},\mathbb{H})\ $for $k=1,2,$\\   
then \ \ 
${|t_kf\left(t\right)|}^2_{2,Q}\ {\left\|{\xi }_k{\mathcal{F}^{\lambda,\mu }}\left\{f\left(t\right)\right\}(2\pi\xi )\right\|}^2_{2,Q}\ \ge \frac{1}{{16\pi }^2}{|f\left(t\right)|}^4_{2,Q}+{COV}^2_{t_k},$ \\
with ${COV}_{t_k}:= \frac{1}{2\pi } \int_{\mathbb R{^2}}{{ {|f\left(t\right)|}}^2_Q{|t_k\left(\frac{\partial }{\partial t_k}e^{u\left(t\right)\theta (t)}\right)|}_Q dt}.$  \\
The equation holds if and only if 
$f\left(t\right)= De^{-a_kt^2_k}e^{u\left(t\right)\theta (t)} $ and  $\frac{\partial }{\partial t_k}e^{u\left(t\right)\theta (t)}= \delta_kt_k,$\\
 where ${a_1,a}_2>0, D\in $  $\mathbb R{^+}$ and $\delta_1, \delta_2\ $are pure quaternions.

\end{Prop}
\begin{Th}\label{QOLCT_Heisenberg}
Suppose that  $f, \ \frac{\partial }{\partial t_k} f,\ t_k f\in L^2({\mathbb R}^2,{\mathbb H})\ $\ for  $k=1,2,$ \\
   then 

${\left|t_kf\left(t\right)\right|}^2_{2,Q}\  {\left\|\frac{{\xi }_k}{2\pi b_k}{\mathcal O}^{\lambda,\mu }_{A_1,A_2}\left\{f\left(t\right)\right\}\left(\xi \right)\right\|}^2_{2,Q}\ge \frac{1}{{16\pi }^2}{\left|f\left(t\right)\right|}^4_{2,Q}+{COV}^2_{t_k\xi},$\\
Where
${COV}_{t_k\xi}:=\ \frac{1}{2\pi }( \int_{{\mathbb R}^2}{{({\left|f(t)\right|}_Q)}^2{\left|t_k\left(\frac{\partial }{\partial t_k}e^{u\left(t\right)\theta (t)}\right)\right|}_Q\ dt}),$\\
and \\
$e^{u\left(t\right)\theta (t)}= \frac{1}{{\left|f(t)\right|}_Q}\frac{1}{\sqrt{\lambda}}e^{-\lambda\frac{1}{b_1}{\xi_1(d}_1\tau_1-b_1\eta_1)+\lambda\frac{d_1}{2b_1}{(\xi }^2_1+{\tau }^2_1)+\lambda\frac{1}{b_1}t_1\tau_1+\lambda\frac{a_1}{2b_1}t^2_1} f(t)\ e^{\mu\frac{1}{b_2}t_2\tau_2+\mu\frac{a_2}{2b_2}t^2_2-\mu\frac{1}{b_2}{\xi_2(d}_2\tau_2-b_2\eta_2)+\mu\frac{d_2}{2b_2}{(\xi }^2_2+{\tau }^2_2)}\frac{1}{\sqrt{\mu}}.$

The equation holds if and only if

$f\left(t\right)= De^{-a_kt^2_k}e^{u\left(t\right)\theta (t)} $ and  $\frac{\partial }{\partial t_k}e^{u\left(t\right)\theta (t)}= \delta_kt_k$\ where ${a_1,a}_2$$>$0, $D\in $  $\mathbb R{^+}$ and $\delta_1, \delta_2\ $are pure quaternions.

\end{Th}

Proof. Let $h(t)$\ be rewritten as  \eqref{function_h}.

 Since$\frac{\partial }{\partial t_k} f,\ t_k f\in L^2({\mathbb R}^2,{\mathbb H})$, and ${\left|h(t)\right|}_Q\ = \frac{1}{2\pi \sqrt{b_1b_2}}{\left|f(t)\right|}_Q.$\\
 we get $\frac{\partial }{\partial t_k}h,\ t_k h\in L^2({\mathbb R}^2,{\mathbb H})$,\ and
${\left|t_k h\left(t\right)\right|}^2_{2,Q} = \int_{{\mathbb R}^2}{{t^2_k  \left|h\left(t\right)\right|}^2_Qdt}= \frac{1}{4{\pi }^2b_1b_2} {\left|t_kf\left(t\right)\right|}^2_{2,Q}$, \\
${\left|h\left(t\right)\right|}^4_{2,Q}=\frac{1}{16{\pi }^4b^2_1b^2_2}{\left|f\left(t\right)\right|}^4_{2,Q}$

By lemma \ref{QFT_QOLCT}, we have

${\left\|{\xi }_k {\mathcal F}^{\lambda,\mu }\left\{h\left(t\right)\right\}\left(2\pi \xi \right)\right\|}^2_{2,Q}={\left\|{\xi }_k {\mathcal O}^{\lambda,\mu }_{A_1,A_2}\left\{f\left(t\right)\right\}\left(2\pi {b_1\xi }_1,{{2\pi b}_2\xi }_2\right)\right\|}^2_{2,Q}$\\
\hspace*{3.5 cm}=$\frac{1}{{4\pi }^2b_1b_2}{\left\|  \frac{{\xi }_k}{2\pi b_k}  {\mathcal O}^{\lambda,\mu }_{A_1,A_2}\left\{f\left(t\right)\right\}\left(\xi \right)\right\|}^2_{2,Q}.$ \\
Hence, it follows from proposition \ref{QFT_Heisenberg}.\\
\begin{equation}\label{eq_Heisenberg}({\left|t_kh\left(t\right)\right|}^2_{2,Q})({\xi}^2_k{\left\|{\mathcal F}^{\lambda,\mu }\left\{h\left(t\right)\right\}(2\pi \xi )\right\|}^2_{2,Q}))\ge \frac{1}{{16\pi }^2}{\left\|h\left(t\right)\right\|}^4_{2,c}+{COV}^2_{t_k},\end{equation}
With  ${COV}_{t_k}~=\frac{1}{2\pi }\ \int_{{\mathbb R}^2}{{\left|h\left(t\right)\right|}^2_Q{\left|t_k\left(\ \frac{\partial }{\partial t_k}e^{u\left(t\right)\theta (t)}\right)\right|}_Q\ dt},$  

and $e^{u\left(t\right)\theta (t)}=\frac{1}{{\left|h(t)\right|}_Q}\ h\left(t\right)$\\
    \hspace*{2.8 cm}=$\frac{2\pi \sqrt{b_1b_2}}{{\left|f(t)\right|}_Q}\frac{1}{\sqrt{2\pi \lambda {b}_1}}e^{-\lambda\frac{1}{b_1}{\xi_1(d}_1\tau_1-b_1\eta_1)+\lambda\frac{d_1}{2b_1}{(\xi}^2_1+{\tau }^2_1)+\lambda\frac{1}{b_1}t_1\tau_1+\lambda\frac{a_1}{2b_1}t^2_1}\ f(t)$\\
	\hspace*{3 cm}$\times e^{\mu\frac{1}{b_2}t_2\tau_2+\mu\frac{a_2}{2b_2}t^2_2-\mu\frac{1}{b_2}{\xi_2(d}_2\tau_2-b_2\eta_2)+\mu\frac{d_2}{2b_2}{(\xi}^2_2+{\tau }^2_2)}\frac{1}{\sqrt{2\pi\mu {b}_2}}$\\
  \hspace*{2.8cm}=$\frac{1}{{\left|f(t)\right|}_Q}\frac{1}{\sqrt{\lambda}}e^{-\lambda\frac{1}{b_1}{\xi_1(d}_1\tau_1-b_1\eta_1)+\lambda\frac{d_1}{2b_1}{(\xi}^2_1+{\tau }^2_1)+\lambda\frac{1}{b_1}t_1\tau_1+\lambda\frac{a_1}{2b_1}t^2_1}\ f(t)$\\ 
  \hspace*{3 cm}$\times e^{\mu\frac{1}{b_2}t_2\tau_2+\mu\frac{a_2}{2b_2}t^2_2-\mu\frac{1}{b_2}{\xi_2(d}_2\tau_2-b_2\eta_2)+\mu\frac{d_2}{2b_2}{(\xi}^2_2+{\tau }^2_2)}\frac{1}{\sqrt{\mu}}.$

\eqref{eq_Heisenberg} implies that\\
${\left|t_k f\left(t\right)\right|}^2_{2,Q} {\left\|\frac{{\xi }_k}{2\pi b_k} {\mathcal O}^{\lambda,\mu }_{A_1,A_2}\left\{f\left(t\right)\right\}\left(\xi\right)\right\|}^2_{2,Q} \ge \frac{1}{{16\pi }^2}{\left|f\left(t\right)\right|}^4_{2,Q}+16{\pi }^4{b^2_1b^2_2\ COV}^2_{t_k}.$ 

After straightforward calculation one obtains\\
$16{\pi }^4{b^2_1b^2_2COV}^2_{t_k}= \frac{1}{{4\pi }^2}{(\ \int_{{\mathbb R}^2}{{({\left|f(t)\right|}_Q)}^2{\left|t_k\left(\frac{\partial }{\partial t_k}e^{u\left(t\right)\theta (t)}\right)\right|}_Q dt})}^2.$\\
By proposition \ref{QFT_Heisenberg}, the equation holds in \eqref{eq_Heisenberg} if and only if $\frac{\partial }{\partial t_k}e^{u\left(t\right)\theta (t)}= \delta_kt_k,$ and 
 ${\left|h\left(t\right)\right|}_Q=Ce^{-a_kt^{{\rm 2}}_k},$\\
 that  is ${\left|f(t)\right|}_Q=2\pi \sqrt{b_1b_2} C e^{-a_kt^{{\rm 2}}_k},$\\
     where  $a_1,a_2 > 0, C,D \in  {\mathbb R}^+$ and $\delta_1,\delta_2\ $are pure quaternions.
This proves the theorem.$\hfill \square$

\subsection{Hardy's uncertainty principle}
Hardy's theorem \cite{HA33} is a qualitative UP, it states that it is impossible for a function and its Fourier transform to decrease rapidly simultaneously.
The following proposition is the Hardy's UP for the Two-sided QFT.
\begin{Prop}\label{Hardy}({\cite{EF17}, Thm. 5.3})\\
Let $\alpha \ $  and $\beta $  are positive constants .Suppose $f\in L^{1}({\mathbb R}^2,{\mathbb H})$  with 

\[{|f\left(t\right)|}_Q\le {Ce}^{-\alpha {\left|t\right|}^2},\  t\in {\mathbb R}^2.\] 

\[{|\ {\mathcal F}^{\lambda,\mu }\left\{f\right\}\left(u\right)|}_Q\le {C' e}^{-\beta {\left|u\right|}^2},\ u\in {\mathbb R}^2.\] 

for some positive constants $C,C'.$Then, three cases can occur :

\end{Prop}

On the basis of proposition \ref{Hardy}, we give the Hardy's UP in the QOLCT domains.
\begin{Th}\label{QOLCT_Hardy}
Let $\alpha  $  and $\beta $  are positive constants. Suppose $f\in L^{1}({\mathbb R}^2,{\mathbb H})$  with 
\begin{equation}\label{Hardy1}{\left|f(t)\right|}_Q\le {Ce}^{-\alpha {\left|t\right|}^2},  t\in {\mathbb R}^2. \end{equation}
\begin{equation}\label{Hardy2}{\left|{\mathcal O}^{\lambda,\mu }_{A_1,A_2}f(b_1u_1,b_2u_2)\right|}_Q \le C'e^{-\beta {\left|u\right|}^2},\ u\in {\mathbb R}^2. \end{equation} 
for some positive constants $C, C'.$ Then, three cases can occur :\\
 
\end{Th}
Proof. Let $g(t)$\ be rewritten in the form of \eqref{function_g}, we have\\
 ${\mathcal O}^{\lambda,\mu }_{A_1,A_2}\left\{f(t)\right\}\left(u_1, u_2\right)=\frac{1}{\sqrt{2\pi {\lambda}_1}}e^{-\lambda\frac{1}{b_1}{u_1(d}_1{\tau}_1-b_1\eta_1)+\lambda\frac{d_1}{2b_1}{(u}^2_1+{\tau}^2_1)}{\mathcal F}^{\lambda,\mu }\left\{g(t)\right\}\left(\frac{u_1}{b_1},\frac{u_2}{b_2}\right)$\\ 

\begin{equation}\label{QFT_QOLCT_g}
\times e^{-\mu\frac{1}{b_2}{u_2(d}_2\tau_2-b_2\eta_2)+\mu\frac{d_2}{2b_2}{(u}^2_2+{\tau }^2_2)}\frac{1}{\sqrt{2\pi \mu{b}_2}}.
\end{equation}
Since
\[\ {\left|g(t)\right|}_Q\ ={\left|f(t)\right|}_Q,\] 
We get  $g \in L^{1}({\mathbb R}^2,{\mathbb H})$   and  $\ {\left|g(t)\right|}_Q\ \le {Ce}^{-\alpha {\left|t\right|}^2}.$\\
On the other hand, by \eqref{QFT_QOLCT_g} and \eqref{Hardy2} we obtain
\[{\left|{\mathcal F}^{\lambda,\mu }\left\{g(t)\right\}\left(u_1{,\ }u_2\right)\right|}_Q=2\pi \sqrt{b_1b_2}{\ \left|{\mathcal O}^{\lambda,\mu }_{A_1,A_2}\left\{f(t)\right\}\left(b_1u_1{,\ }b_2u_2\right)\right|}_Q\] 
\[\le 2\pi \sqrt{b_1b_2}C'e^{-\beta {\left|u\right|}^2}.\] 
Therefore, it follows from proposition \ref{Hardy} that,\\ 
If $\alpha \beta =\frac{1}{4}$ then\\
$g\left(t\right)=Ae^{-\alpha {\left|t\right|}^2}$, for some  constant $A.$\\
Hence 
\[f\left(t\right)=A e^{- \alpha {\left|t\right|}^2}  e^{-\lambda\frac{1}{b_1}t_1\tau_1-\lambda\frac{a_1}{2b_1}t^2_1}e^{-\mu\frac{1}{b_2}t_2\tau_2-\mu\frac{a_2}{2b_2}t^2_2}.\] 
If $\alpha \beta >\frac{1}{4}$ then $g=0$, so $f=0$.\\
If $\alpha \beta <\frac{1}{4},$ then there are infinitely many such functions $f$, that verify  \eqref{Hardy1} and ( \eqref{Hardy2}.\\
This completes the proof.$\hfill \square$\\
It follows from theorem \ref{QOLCT_Hardy} that it is impossible for $f$\ and its two-sided QOLCT to both decrease very rapidly.

\subsection{Beurling's uncertainty principle}

Beurling's UP \cite{BE89}, \cite{HO91} is a variant of Hardy's UP. It implies the weak form of Hardy's UP immediatly.
The following proposition is the Beurling's UP for the Two-sided QFT.
\begin{Prop}\label{Beurling}{ \cite{FE17}}\\
 
Let $f \in L^2\left(\mathbb {\mathbb R}^2,\mathbb H\right)\ and\ \ d\ge 0~$   satisfy  

\hspace*{4 cm} $\int_{\mathbb {\mathbb R}^2}{\int_{\mathbb {\mathbb R}^2}{\ \frac{{\left|f\right|}_Q\  {\left\|\mathcal F\left\{f\right\}\left(y\right)\right\|}_Q}{{(1+\left|x\right|+\left|y\right|)}^d}e^{2\pi \left|x\right|\left|y\right|}}}\ dxdy <\infty, $  

Then $f(x)= P(x)e^{-a{\left|x\right|}^2},$\ \ a.e. 

Where $a>0$ and $P$ is a polynomial of degree $<\frac{d-2}{2}$. In particular, $f$ is identically 0  when $d\le 2.$

\end{Prop}
On the basis of proposition \ref{Beurling}, we give the Beurlings'UP in the QOLCT domains.

\begin{Th}
Let $f \in L^2\left(\mathbb {\mathbb R}^2,\mathbb H\right)\ and\ \ d\ge 0~$   satisfy  

$\int_{{\mathbb R}^2}{\int_{{\mathbb R}^2}{ \frac{{\left| f(t )\right|}_Q\ {\left\|{\mathcal O}^{\lambda,\mu }_{A_1,A_2}\left\{f\right\}\left(b_1u_1,b_2u_2\right)\right\|}_Q}{{(1+\left|t\right|+\left|u\right|)}^d}e^{\left|t\right|\left|y\right|}}}\ dtdu <\infty ,$ \\   
Then \\ 
$f\left(t\right)=e^{-a{\left|t\right|}^2}{\sqrt{2\pi {\lambda b}_1}e}^{\lambda\frac{1}{b_1}{u_1(d}_1\tau_1-b_1{\eta}_1)-\lambda\frac{d_1}{2b_1}{(u}^2_1+{\tau }^2_1)-\lambda\frac{1}{b_1}t_1\tau_1-\lambda\frac{a_1}{2b_1}t^2_1} P(t)\ e^{\mu\frac{1}{b_2}{u_2(d}_2\tau_2-b_2\eta_2)-\mu\frac{d_2}{2b_2}{(u}^2_2+{\tau }^2_2)-\mu\frac{1}{b_2}t_2\tau_2-\mu\frac{a_2}{2b_2}t^2_2}$\\ 
\hspace*{6cm}$\times\sqrt{2\pi\mu {b}_2},$ \ \ a.e. 

Where $a>0$ and $P$\ is a quaternion polynomial of degree $< \frac{d-2}{2}$. In particular, $f=0$\  a.e. when $d\le 2.$
\end{Th}

Proof. Let $h(t)$\ be rewritten in the form of  \eqref{function_h},\ we have $h \in L^2\left(\mathbb {\mathbb R}^2,\mathbb H\right).$\\ 
It follows from lemma \ref{QFT_QOLCT}, and ${\left|h\left(t\right)\right|}_Q=\frac{1}{2\pi \sqrt{b_1b_2}}{\left|\ f(t)\ \right|}_Q$
That, \\
\hspace*{1 cm}$\int_{{\mathbb R}^2}{\int_{{\mathbb R}^2}{ \frac{{\left|h(t)\right|}_Q {\left\|{\mathcal F}^{\lambda,\mu }\left\{h\right\}\left(u\right)\right\|}_Q}{{(1+\left|t\right|+\left|u\right|)}^d}e^{\left|t\right|\left|u\right|}}}\ dtdu=\int_{{\mathbb R}^2}{\int_{{\mathbb R}^2}{ \frac{{\left|h(t)\right|}_Q {\left\|{\mathcal O}^{\lambda,\mu }_{A_1,A_2}\left\{f\right\}\left(b_1u_1,b_2u_2\right)\right\|}_Q}{{(1+\left|t\right|+\left|u\right|)}^d}e^{\left|t\right|\left|u\right|}}} dtdu$\\
\hspace*{6.7 cm}=$\frac{1}{2\pi \sqrt{b_1b_2}}\int_{{\mathbb R}^2}{\int_{{\mathbb R}^2}{ \frac{{\left|f(t)\right|}_Q {\left\|{\mathcal O}^{\lambda,\mu }_{A_1,A_2}\left\{f\right\}\left(b_1u_1,b_2u_2\right)\right\|}_Q}{{(1+\left|t\right|+\left|u\right|)}^d}{\ e}^{\left|t\right|\left|u\right|}}}\ dtdu<\infty .$\\
Then by proposition \ref{Beurling}, we get
 $h(t)= P(t) e^{-a{\left|t\right|}^2}$ a.e.
where $a>0$\ and $P$\ is a quaternion polynomial of degree $<\frac{d-2}{2}.$\\
i.e.\\
$f\left(t\right)=e^{-a{\left|t\right|}^2}{\sqrt{2\pi \lambda{b}_1}e}^{\lambda\frac{1}{b_1}{u_1(d}_1\tau_1-b_1\eta_1)-\lambda\frac{d_1}{2b_1}{(u}^2_1+{\tau }^2_1)-\lambda\frac{1}{b_1}t_1\tau_1-\lambda\frac{a_1}{2b_1}t^2_1} P(t) e^{\mu\frac{1}{b_2}{u_2(d}_2\tau_2-b_2\eta_2)-\mu\frac{d_2}{2b_2}{(u}^2_2+{\tau }^2_2)-\mu\frac{1}{b_2}t_2\tau_2-\mu\frac{a_2}{2b_2}t^2_2}$\\
\hspace*{6.7 cm}$\times\sqrt{2\pi \lambda{b}_2}.$\\
In particular, $f=0$\ a.e. when $d \le 2.\hfill \square$
\subsection{Logarithmic uncertainty principle}

The logarithmic UP \cite{BE95} is a more general form of Heisenberg type UP, its localization is measured in terms of entropy. It is derived by using Pitt's inequality.
\begin{Lemma}\label{Pitt}{ (Pitt's inequality for the two-sided QFT \cite{CKL15})}\\
 For  $f\in {\mathcal S}({\mathbb R}^2,{\mathbb H})$, and   $0\le \alpha <2$,
\[\int_{{\mathbb R}^2}{ {\left|u\right|}^{-\alpha }}{\left\|{\mathcal F}^{i,j}\left\{f(t)\right\}\left(u_1,u_2\right)\right\|}^2_Qdu\  \le C_{\alpha } \int_{{\mathbb R}^2}{ {\left|t\right|}^{\alpha }}{\left|f(t)\right|}^2_Q\ dt.\] 
With $C_{\alpha }:=\frac{{4\pi }^2}{2^{\alpha }}{{ [}\Gamma (\frac{2-\alpha }{4})/\Gamma (\frac{2+\alpha }4)]}^2$, and $\Gamma \left(.\right)$\ is the Gamma function.
\end{Lemma}

\begin{Th}
Under the assumptions of lemma \ref{Pitt}, one has
\begin{equation}\label{Becker1}\int_{{\mathbb R}^2}{ {\left|(\frac{z_1}{b_1},\frac{z_2}{b_2})\right|}^{-\alpha }}{\left\|{\mathcal O}^{i,j}_{A_1,A_2}\{f \}\left(z_1,z_2\right)\right\|}^2_Q dz\  \le \ \frac{C_{\alpha }}{{4\pi }^2}\ \int_{{\mathbb R}^2}{ {\left|t\right|}^{\alpha }}{\left|f(t)\right|}^2_Q\ dt. \end{equation} 
\end{Th}

Proof. Let $h(t)$ be rewritten in the form of  \eqref{function_h}, with $\lambda=i,$\ and $\mu=j$.\\ 
It's clear that  $h \in {\mathcal S}\left(\mathbb {\mathbb R}^2,\mathbb H\right)$, and  $\ {\left|h(t)\right|}_Q = \frac{1}{2\pi \sqrt{b_1b_2}}{\left|f(t)\right|}_Q.$ \\
Let  ${\mathcal O}^{i,j}_{A_1,A_2}\left\{f\right\}\left(u\right)$ be rewritten as \eqref{qolct}, 
we have by lemma \ref{QFT_QOLCT}

\[{\mathcal O}^{i,j}_{A_1,A_2}\left\{f\left(t\right)\right\}\left(u_1,u_2\right)= {\mathcal F}^{i,j}\left\{h(t)\right\}\left(\frac{u_1}{b_1},\frac{u_2}{b_2}\right). \] 

By lemma \ref{Pitt}, we obtain 
\[\int_{{\mathbb R}^2}{ {\left|u\right|}^{-\alpha }}{\left\|{\mathcal O}^{i,j}_{A_1,A_2}\left\{f\right\}\left({b_1u}_1{,\ }{b_2u}_2\right)\right\|}^2_Qdu = \int_{{\mathbb R}^2}{\ \ {\left|u\right|}^{-\alpha }}{\left\|{\mathcal F}^{i,j}\left\{h(t)\right\}\left(u_1,u_2\right)\right\|}^2_Qdu\] 
\hspace*{9.5cm}$\le  C_{\alpha }\int_{{\mathbb R}^2}{ {\left|t\right|}^{\alpha }}{\left|h(t)\right|}^2_Q\ dt= \frac{C_{\alpha }}{{4\pi }^2b_1b_2}\ \int_{{\mathbb R}^2}{ {\left|t\right|}^{\alpha }}{\left|f(t)\right|}^2_Q\ dt.$ \\
Let\ $z_1=b_1u_1$\ and $z_2=b_2u_2$, we have
\[\frac{1}{b_1b_2}\int_{{\mathbb R}^2}{ {\left|(\frac{z_1}{b_1},\frac{z_2}{b_2})\right|}^{-\alpha }}{\left\|{\mathcal O}^{i,j}_{A_1,A_2}\left\{f\right\}\left(z\right)\right\|}^2_Qdz\le \ \frac{C_{\alpha }}{{{4\pi }^2b}_1b_2}\ \int_{{\mathbb R}^2}{ {\left|t\right|}^{\alpha }}{\left|f(t)\right|}^2_Q dt,\] 
i.e.,
\[\int_{{\mathbb R}^2}{ {\left|(\frac{z_1}{b_1},\frac{z_2}{b_2})\right|}^{-\alpha }}{\left\|{\mathcal O}^{i,j}_{A_1,A_2}\left\{f\right\}\left(z\right)\right\|}^2_Qdz\le \ \frac{C_{\alpha }}{{4\pi }^2}\ \int_{{\mathbb R}^2}{ {\left|t\right|}^{\alpha }}{\left|f(t)\right|}^2_Q dt.\] $\hfill \square$

\begin{Th}{(Logarithmic UP for the QOLCT)}\\
Let $f \in {\mathcal S}\left({\mathbb R}^2,{\mathbb H}\right), $\ then
\begin{equation}\label{Becker2}\int_{{\mathbb R}^2}{\ {{ ln(} \left|\frac{z_1}{b_1},\frac{z_2}{b_2}\right|)\ \ }}{\left\|{\mathcal O}^{i,j}_{A_1,A_2}\left\{f\right\}\left(z\right)\right\|}^2_Qdz+\int_{{\mathbb R}^2}{\ {\ln  \left(\left|t\right|\right)\ }\ }{\left|f(t)\right|}^2_Q\ dt\ge \  A\int_{{\mathbb R}^2}{\ \ }{\left|f(t)\right|}^2_Q\ dt,\end{equation}
with   $A={\ln  \left(2\right)}+{\Gamma '\left(\frac{1}{2}\right)}/{\Gamma (\frac{1}{2})}.$

\end{Th}
Proof. Let $f \in {\mathcal S}\left({\mathbb R}^2,{\mathbb H}\right)$, \   $0\le \alpha <4$,
 $D_{\alpha }=\frac{C_{\alpha }}{{4\pi }^2}\ =\frac{1}{2^{\alpha \ }}{{ [}\Gamma (\frac{2-\alpha }{4})/\Gamma (\frac{2+\alpha }{4})]}^2,$\\
and
$\Phi (\alpha ) := \int_{{\mathbb R}^2}{ {\left|(\frac{z_1}{b_1},\frac{z_2}{b_2})\right|}^{-\alpha }}{\left\|{\mathcal O}^{i,j}_{A_1,A_2}\left\{f\right\}\left(z\right)\right\|}^2_Qdz-D_{\alpha }\ \int_{{\mathbb R}^2}{{\left|t\right|}^{\alpha }}{\left|f(t)\right|}^2_Q\ dt.$\\
By differentiating  $\Phi (\alpha)$, we have

$\Phi '(\alpha )=-\int_{{\mathbb R}^2}{ {{ ln(} \left|\frac{z_1}{b_1},\frac{z_2}{b_2}\right|)}{\left|(\frac{z_1}{b_1},\frac{z_2}{b_2})\right|}^{-\alpha }}{\left\|{\mathcal O}^{i,j}_{A_1,A_2}\left\{f\right\}\left(z\right)\right\|}^2_Qdz-{ D}_{\alpha }'\int_{{\mathbb R}^2}{ {\left|t\right|}^{\alpha }}{\left|f(t)\right|}^2_Qdt-D_{\alpha }\int_{{\mathbb R}^2}{{\ln  \left(\left|t\right|\right)}{\left|t\right|}^{\alpha }}{\left|f(t)\right|}^2_Q\ dt,$\\
whith\\
 ${D}_{\alpha }'=- \ln(2) 2^{-\alpha }{[\Gamma (\frac{2-\alpha }{4})/\Gamma (\frac{2+\alpha }{4}){ ]}}^2+2^{-\alpha } {[-\frac{1}{2}\Gamma (\frac{2-\alpha }{4}){\Gamma '}\left(\frac{2-\alpha }{4}\right){\Gamma }^2\left(\frac{2+\alpha }{4}\right)-\frac{1}{2}{\Gamma }^2\left(\frac{2-\alpha }{4}\right)\Gamma (\frac{2+\alpha }{4})\Gamma { '}\left(\frac{2+\alpha }{4}\right)]}\ /{{\Gamma }^4(\frac{2+\alpha }{4})}.$\\
We have ${\ D}_0=1 $\ and\ $D_0'=-{\ln \left(2\right) }-{\Gamma '\left(\frac{1}{2}\right)}/{\Gamma (\frac{1}{2})}.$\\
Because of \eqref{Becker1}, we see that
$\Phi (\alpha )\le 0$\ for  $0\le \alpha <2$,
also by theorem \ref{Plancherel_QOLCT} we have $\Phi \left(0\right)=0.$\\
Then ${\Phi}^{'}\left(0^+\right)={\mathop{\lim }_{\alpha \to 0^+} \frac{\Phi \left(\alpha \right)-\Phi (0)}{\alpha }}\le 0.$\\
Therefore
 $({\ln  \left(2\right)\ }+{\Gamma'\left(\frac{1}{2}\right)}/{\Gamma (\frac{1}{2})}$)$\ \int_{{\mathbb R}^2}{\left|f(t)\right|}^2_Q\ dt
\le \int_{{\mathbb R}^2}{ {ln( \left|\frac{z_1}{b_1},\frac{z_2}{b_2}\right|) }}{\left\|{\mathcal O}^{i,j}_{A_1,A_2}\left\{f\right\}\left(z\right)\right\|}^2_Qdz+\int_{{\mathbb R}^2}{\ {\ln  \left(\left|t\right|\right) } }{\left|f(t)\right|}^2_Q\ dt.$
From which the theorem follows.$\hfill \square$

\begin{Rem}
Applying Jensen's inequality to \eqref{Becker2}, we can show that the logarithmic UP implies Heisenberg-Weyl's inequality (theorem \ref{QOLCT_Heisenberg}) .
\end{Rem}

\section{Conclusion}
In this paper, we first  presented a new generalization of the QLCT and so of the QFT, namely the QOLCT. Second, We established some properties of the QOLCT including the Plancherel's formula. Then, we derive three UPs in the QOLCT domain: Heisenberg-Weyl's UP, Hardy's UP and its variant-Beurling's UP.
These three UPs assert that it is impossible for a non-zero function and its QOLCT to both decrease very rapidly. Finally, we generalize Pitt's inequality to the QOLCT domain, and then obtain a logarithmic UP  associated with QOLCT.
In the future work, we will consider these UPs for the offset linear canonical transform in Clifford analysis.


\newpage
Let $f$\ and $g\in L^1({\mathbb R}^2,{\mathbb H})$,\ the constants $\alpha $\ and $\beta \in {\mathbb R}, u \in {\mathbb R}^2,$\\$ A_l=\left[\left| \begin{array}{cc}
a_l & b_l \\ 
c_l & d_l \end{array}
\right| \begin{array}{c}
{\tau }_l \\ 
{\eta }_l \end{array}
\right]$, parameters  $a_l,b_l,c_l,d_l,{\tau }_l,{\eta }_l\in {\mathbb R}$ such that $a_ld_l-b_lc_l=1$, for $l=1,2.$\\

\begin{table}[htb]
\begin{tabular}{|p{1.4in}|p{1.2in}|p{2.3in}|} \hline 
\multicolumn{3}{|p{4.5in}|}{ Property   \hspace*{2.3cm}                    Function       \hspace*{2.5cm}            QOLCT} \\ \hline 
Real linearity  \newline  Plancherel's identity & $\alpha f+\beta g$      \newline           ${\left|f\right|}_{Q,2}=$\newline  & $\alpha{\mathcal O}^{\lambda ,\mu }_{A_1,A_2}\left\{f\right\}\left(u\right){\rm +}\beta {\mathcal O}^{\lambda ,\mu}_{A_1,A_2}\left\{g\right\}(u)$ \newline ${\left\|{\mathcal O}^{\lambda,\mu }_{A_1,A_2}\left\{f\right\}\right\|}_{Q,2}$ \\ \hline 
\end{tabular}
\caption{\label{QOLCT_prop} Properties of the quaternionic offset linear canonical transform (QOLCT).}
\end{table}

\end{document}